\documentclass{article}
\usepackage[final]{arxiv}

\PassOptionsToPackage{numbers, compress}{natbib}
\usepackage{natbib}

\usepackage{booktabs} 
\usepackage{algorithm}
\usepackage{algpseudocode}
\usepackage{todonotes} 
\usepackage{placeins} 
\usepackage{amsmath}
\usepackage{amssymb}
\usepackage{epstopdf}
\usepackage{wrapfig}
\usepackage{subfigure}
\usepackage{xspace}
\usepackage{pythonhighlight}
\usepackage[hidelinks]{hyperref}

\newcommand{\bbR}{\mathbb{R}}

\newcommand{\bbZ}{\mathbb{Z}}
\newcommand{\calD}{\mathcal{D}}
\newcommand{\calF}{\mathcal{F}}
\newcommand{\calM}{\mathcal{M}}
\newcommand{\calO}{\mathcal{O}}
\newcommand{\MATLAB}{\textsc{Matlab}}
\newcommand{\CPP}{C\texttt{++}}

\algnewcommand{\algorithmicand}{\textbf{ and }}
\algnewcommand{\algorithmicor}{\textbf{ or }}
\algnewcommand{\OR}{\algorithmicor}
\algnewcommand{\AND}{\algorithmicand}
\algnewcommand{\algorithmicgoto}{\textbf{go to}}%
\algnewcommand{\Goto}[1]{\algorithmicgoto~\ref{#1}}%

\newcommand{\pysot}{\texttt{pySOT}\xspace}
\newcommand{\poap}{\texttt{POAP}\xspace}

\graphicspath{{./pics/}}
\begin{document}

\title{{py}SOT and POAP: An event-driven asynchronous framework for surrogate optimization}

\author{
  David Eriksson, Cornell University$^*$ \\ 
  David Bindel, Cornell University \\
  Christine A. Shoemaker, National University of Singapore
}

{\let\thefootnote\relax\footnote{
  {$^*$The work was conducted while David Eriksson was at Cornell University --- he is currently a Research Scientist at Uber AI Labs}}}

\maketitle

\begin{abstract}
%
%
This paper describes Plumbing for Optimization with Asynchronous Parallelism
(\poap) and the Python Surrogate Optimization Toolbox (\pysot).
%
%
\poap is an event-driven framework for building and combining asynchronous
optimization strategies, designed for global optimization of expensive
functions where concurrent function evaluations are useful.
\poap consists of three components: a worker pool
capable of function evaluations, strategies to propose evaluations or other
actions, and a controller that mediates the interaction between the workers and strategies.
%
%
\pysot is a collection of synchronous and asynchronous surrogate optimization
strategies, implemented in the \poap framework. We support the stochastic
RBF method \citet{regis2007stochastic} along with various extensions of this method,
and a general surrogate optimization strategy that covers most Bayesian optimization methods.
We have implemented many different surrogate models, experimental designs, acquisition functions,
and a large set of test problems.
%
%
We make an extensive comparison between synchronous and asynchronous parallelism and find that the advantage of 
asynchronous computation increases as the variance of the evaluation time or number of processors increases. 
We observe a close to linear speed-up with 4, 8, and 16 processors in both the synchronous and asynchronous setting.
\end{abstract}

\section{Introduction}
\label{sec:intro}
\subsection{Problem Statement}

%
%

We consider the global optimization problem
\begin{equation}
    \label{eq:globprob}
    \begin{array}{lll}
		\text{minimize}   & f(x)  		   				  								\\
		\text{subject to} & x \in \mathcal{D} \cap (\bbZ^q \times \bbR^{d-q})
    \end{array}
\end{equation}
where $f: \bbZ^q \times \bbR^{d-q} \to \bbR$ is a computationally
expensive black-box function. We assume in addition
that $\mathcal{D}$ is a compact hypercube and that $f(x)$ is a
continuous function over the continuous variables. In our setting, $f(x)$ is non-linear, and has multiple local minima,
and the gradient of $f(x)$ is not available. Computationally expensive refers to any problem where a
single function evaluation takes anywhere between a few minutes and many hours. Common examples include
running an expensive simulation model of a complex physical process and tuning machine learning models \cite{snoek2012practical}.
It is common to have limited time and evaluation budgets due to the significant amount of time necessary for
each function evaluation, making it challenging to find a good solution to (\ref{eq:globprob}) in the
case when $f$ is multimodal.

\subsection{Survey of Methods}

%
%

Many popular algorithms for black-box optimization are not suitable when the
function evaluations are computationally expensive.
Derivative based methods are appealing in cases when gradient information
can be obtained cheaply, in which case it is possible to run a local optimizer
with a multi-start strategy such as Newton's method or BFGS \cite{avriel2003nonlinear}.
Finite differences can be used when gradient information
is unavailable, but it is very computationally expensive since $f(x)$ is expensive, and impreciseness in the simulation model often leads to
inaccurate estimates. Several popular derivative free optimization (DFO) methods exist for local optimization such as pattern search \cite{hooke1961direct}, Nelder-Mead \cite{nelder1965simplex}, and ORBIT \cite{wild2011global},  but these methods are not good choices
for multimodal problems. Global heuristic optimization methods such as genetic algorithms \cite{goldberg2006genetic},
particle swarm optimization \cite{kennedy2010particle}, and differential
evolution \cite{storn1997differential}, generally
require a large number of function evaluations and
are not practical for computationally expensive objective functions.

%
%

A successful family of optimization algorithms for computationally expensive
optimization are methods based on surrogate models. The surrogate model approximates
the objective function and helps accelerate convergence to a good solution.
Popular choices are methods based on radial basis functions (RBFs) such as
\cite{regis2007stochastic,regis2013combining,gutmann2001radial}
and Kriging and Gaussian process (GPs) based methods such as
\cite{jones2001taxonomy,jones1998efficient,frazier2008knowledge}.
Other possible surrogate models are polynomial
regression models and multivariate adaptive regression splines \cite{friedman1991multivariate,muller2014influence}.
Most surrogate optimization algorithms
start by evaluating an experimental design that is used to fit
the initial surrogate model. What follows is an adaptive phase where an auxiliary
problem is solved to pick the next sample point(s), and this phase continues until
either a restart or a stopping criterion has been met. We can avoid getting trapped in
a local minimum by using an auxiliary problem that provides a good balance of exploration
and exploitation.

%
%

Several parallel algorithms have been developed for computationally expensive
black-box optimization. Regis and Shoemaker \cite{regis2009parallel} developed
a synchronous parallel surrogate optimization algorithm based on RBFs and this
idea was later extended to SOP algorithm for large number of processors \cite{krityakierne2016sop}.
In both algorithms, it is assumed that (i)  the resources are homogeneous and (ii) the evaluation time is constant.
The first assumption does not hold for heterogeneous parallel computing platforms and the second assumption is unlikely
to hold in cases where the complexity of evaluating the objective depends spatially on the input. The first assumption
can almost always be assessed before the start of the optimization run while the second assumption may not be easy
to assess in practice.

\begin{wrapfigure}{R}{0.5\textwidth}
    \includegraphics[width=0.48\textwidth]{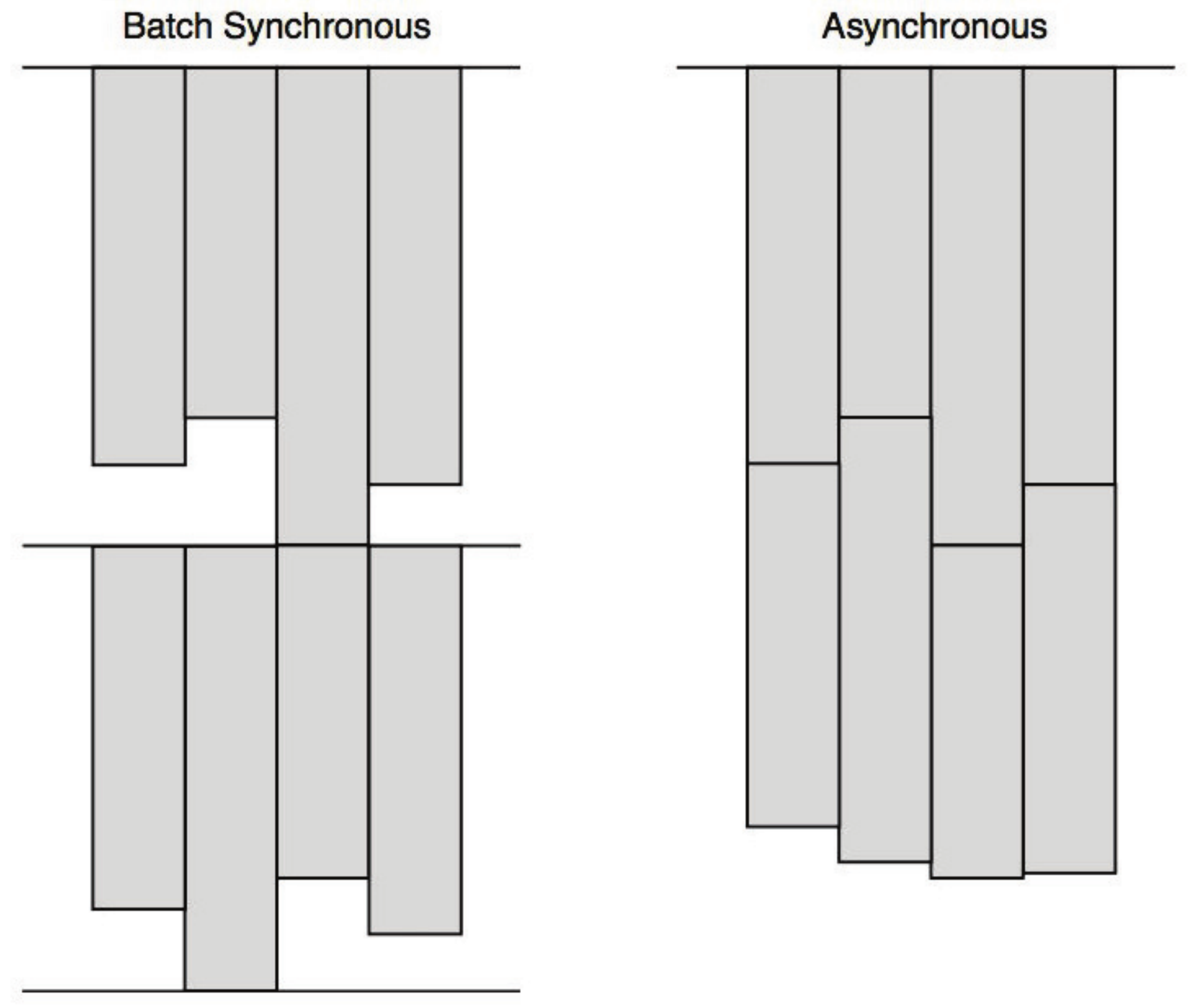}
	\caption{\textit{Synchronous vs asynchronous parallel.}}
	\label{fig:async_vs_sync}
\end{wrapfigure}

Another limitation of the work in \cite{regis2009parallel} is that the algorithm
does not handle the possibility of worker failures and crashed evaluations. Being able
to handle failures is critical in order to run the algorithm on large-scale systems.
The natural
way of dealing with cases where (i) or (ii) are violated is to launch function
evaluations asynchronously, which is illustrated in Figure
\ref{fig:async_vs_sync} to eliminate idle time.

\subsection{Survey of Software}

%
%

A library with similar functionality to \poap is \texttt{SCOOP} \cite{hold2014once},
a Python based library for distributing concurrent tasks while internally
handling the communication. \poap provides similar functionality for global
optimization problems and also handles all of the communication internally,
which makes it easy to implement asynchronous optimization algorithms.

%
%

\texttt{HOPSPACK} (Hybrid Optimization Parallel Search PACKage) \cite{plantenga2009hopspack} is
a \CPP{} framework for derivative-free optimization problems. \texttt{HOPSPACK} supports
parallelism through MPI or multi-threading and supports running multiple optimization
solvers simultaneously, a functionality similar to combining strategies in \poap.
The framework implements an asynchronous pattern search solver and supports non-linear
constraints and mixed-integer variables, but there is no support for surrogate optimization.


%
%

\texttt{MATSuMoTo} (\MATLAB{} Surrogate Model Toolbox) \cite{mueller2014matsumoto}
is an example of a surrogate global optimization
toolbox. \texttt{MATSuMoTo} is written in \MATLAB{} and has support for computationally expensive,
black-box global optimization problems that may have continuous, mixed-integer, or
pure integer variables. \texttt{MATSuMoTo} offers a variety of choices for surrogate models and
surrogate model mixtures, experimental designs, and auxiliary functions.
The framework is not designed to support a large class of surrogate optimization algorithms
and the lack of object orientation makes it hard to extend the framework. Parallelism is
only supported through \MATLAB{}'s Parallel Computing Toolbox and there is no support
for asynchrony, combining strategies, or dynamically changing the number of workers. Furthermore,
many large-scale systems do not support \MATLAB{}. Note that as of version 2018b, the 
\MATLAB{} Global Optimization Toolbox offers \texttt{surrogateopt}
\footnote{\url{https://www.mathworks.com/help/gads/surrogateopt.html}}, which is an 
asynchronous surrogate optimization method implementation based on \citet{regis2007stochastic}.

%
%

Nonlinear Optimization by Mesh Adaptive Direct Search (\texttt{NOMAD})
\cite{le2011algorithm} is a library intended for time-consuming black-box
simulation with a small number of variables. The library implements mesh adaptive direct
search (MADS) and there is support for asynchronous function evaluations using
MPI. The framework is fault resilient in the sense that it supports
objective function failing to return a valid output. Similar fault
resilience is provided by \poap, which allows the user to decide what action to take
in case of a failure.

%
%

\texttt{Dakota} \cite{eldred2007dakota} is an extensive toolkit with algorithms for optimization
with and without gradient information; uncertainty quantification,
nonlinear least squares methods, and sensitivity/variance analysis. These
components can be used on their own or with strategies such as surrogate-based
optimization, mixed integer nonlinear programming, or optimization under
uncertainty. The \texttt{Dakota} toolkit is object-oriented and written in \CPP{} with the
intention of being a flexible and extensible interface between simulation codes,
and there is support for parallel function evaluations. \texttt{Dakota} includes \CPP{} code
for global optimization with a GP based surrogate (e.g, an implementation of
the GP-EI/EGO method \cite{jones1998efficient} and EGRA method \cite{bichon2008efficient}).
\texttt{Dakota} does not have global optimization codes designed to be used with RBF surrogates, although
it is possible to construct an RBF surrogate in \texttt{Dakota}.

%
%

\texttt{BayesOpt} \cite{martinez2014bayesopt} is a library with Bayesian optimization
methods to solve nonlinear optimization problems. Bayesian optimization methods
build a posterior distribution to capture the evidence and prior knowledge of
the target function. Built in \CPP{}, the library is efficient, portable, and flexible.
There is support for commonly used methods such as sequential Kriging optimization
(SKO), sequential model-based optimization (SMBO), and efficient global
optimization (EGO). The software is sequential and there is no support for parallel
function evaluations.

%
%

\texttt{RBFOpt} \cite{costa2014rbfopt} is a radial basis function based library that
implements and extends the global optimization RBF algorithm proposed by Gutmann
\cite{gutmann2001radial}. \texttt{RBFOpt} is written in Python and supports asynchronous
parallelism through Python's multiprocessing library, but there is no support
for MPI. The software is not designed to cover a large class of surrogate
optimization methods and there is no support for dynamically changing the number
of workers and combining different optimization strategies.

%
%

\texttt{Cornell-MOE} is a Python library that implements Bayesian optimization with the
expected improvement and knowledge gradient acquisition functions. The software
is built on work that extends these acquisition functions to batch synchronous
parallel, both with and without gradient information \cite{wu2016parallel,wu2017bayesian}.
There is no support for asynchronous parallelism and it is not possible to dynamically
change the number of workers.

\subsection{Contribution}

%
%

The \poap and \pysot software has become very popular with \pysot having been downloaded more than
88,000 times and \poap downloaded more than 126,000 times. The main contribution of \poap is an
event-driven framework for building and
combining asynchronous optimization strategies. The user can implement their own
strategies that specify what actions to take when different events occur, while
all communication and dispatching of work is handled by the framework. \poap is designed
to be both flexible and easily extensible, and the framework makes it easy to
dynamically change the number of workers and combine different optimization strategies.
\poap is fault resilient and handles function evaluation and worker crashes.

%
%

\pysot is a great test-suite for doing head-to-head
comparisons with different experimental designs, surrogate models, and
acquisition functions. Being built on top of \poap, \pysot leverages the many benefits of
the \poap framework, leading to a robust and flexible framework without having to
worry about the communication and dispatching of work. The object-oriented design
makes \pysot easy to extend and users can experiment with different surrogate models,
experimental designs, and auxiliary problems, and make comparisons in either a
synchronous or an asynchronous setting. In addition, \pysot supports checkpointing
which allows users to resume a crashed optimization run.

%
%

We provide an extensive comparison of synchrony and asynchrony in cases where
the objective function evaluation time varies and conclude that reducing idle time is
more important than information for multimodal problems. We conclude that asynchrony
should be preferred over synchrony in this case. The performance difference between
asynchrony and synchrony increases with function evaluation variance and
number of processors since both increase the idle time for synchrony. Our numerical
experiments also indicate that parallelism improves exploration and that the
parallel algorithms often outperform the serial version with respect to number of
function evaluations.

\subsection{Overview}
We review the general surrogate optimization algorithm and
the most common surrogate models, experimental designs, and auxiliary problems
in \S\ref{sec:background}. We describe in detail our asynchronous surrogate
optimization algorithm in \S\ref{sec:async}. The implementation details of \poap
and \pysot are described in \S\ref{sec:poap} and \S\ref{sec:pysot} respectively.
We illustrate a code example in \S\ref{sec:code} that shows how to use \pysot and \poap.
We provide an extensive comparison between asynchrony and synchrony in
\S\ref{sec:experiments} and conclude in \S\ref{sec:conclusions}.

\section{Surrogate optimization}
\label{sec:background}
Most surrogate optimization methods follow the same main steps, both when
running in serial and synchronous parallel. The first step consists of
generating an experimental design that is evaluated to fit an
initial surrogate model. Once an initial surrogate model has been built,
we proceed to optimize an acquisition function to find new point(s) to evaluate.
We often refer to optimizing this acquisition function as solving an auxiliary problem.
We evaluate the new point(s), update the surrogate model, and repeat this
procedure until a stopping criterion has been met. This is illustrated
in Algorithm \ref{alg:sync_algo}.

\begin{algorithm}
	\caption{Synchronous surrogate optimization algorithm}
	\label{alg:sync_algo}
  	\begin{algorithmic}[1]
        \State Generate an experimental design
        \State Evaluate $f(x)$ at the points in the experimental design
        \State Build a surrogate model of $f(x)$ from the data
        \While{Stopping criterion not met}
            \State Solve an auxiliary problem to select the next point(s) to evaluate
            \State Evaluate the new point(s)
            \State Update the surrogate model
        \EndWhile
  	\end{algorithmic}
\end{algorithm}
The overhead of fitting the surrogate model and optimizing the acquisition function should be compared to the evaluation
time of the objective function, as we want to spend most of the computational time on function evaluations.
We proceed with a brief background to the most popular experimental
designs, surrogate models, and auxiliary problems.

\subsection{Experimental design}
\label{sec:expdes}
The simplest experimental design is choosing the $2^d$ corners of the hypercube $\calD$, often referred to as the
2-factorial design, but this is infeasible when $d$ is large and the function is expensive. Two common
alternatives are the Latin hypercube design (LHD) and the symmetric Latin hypercube design (SLHD), which allow an
arbitrary number of design points. We deal with integer variables by rounding the generated design and generate a
new experimental design if the resulting design is rank-deficient or if any two points coincide. This works well
in practice, which is also reported in \cite{costa2014rbfopt} and \cite{muller2013so}.

\subsection{Surrogate models}
\label{sec:surr}
The surrogate model is used to approximate the objective function.
The surrogate model of choice in \pysot is radial basis functions (RBFs),
but we also support Gaussian processes (GPs),
support vector regression (SVR), multivariate adaptive regression splines
(MARS), and polynomial regression.

\subsubsection{Radial basis functions}
\label{sec:rbf}
RBF interpolation is one of the most popular
approaches for approximating scattered data in a general number of dimensions
\cite{buhmann2003radial,fasshauer2007meshfree,schaback2006kernel,wendland2004scattered}.
Given a set of pairwise distinct interpolation points $X=\{x_i\}_{i=1}^n
\subset \Omega$ the RBF model takes the form
\begin{equation}
  \label{eq:rbf}
  s_{f,X}(x) = \sum_{i=1}^n \lambda_i \varphi(\|x - x_i\|) + p(x)
\end{equation}
where the kernel $\varphi : \mathbb{R}_{\geq 0} \to \mathbb{R}$ is a
one-dimensional function and $p \in \Pi_{k-1}^d$, the space of
polynomials with $d$ variables of degree no more than $k-1$.
The name RBF comes from the fact that the function
$\varphi(\cdot)$ is constant on spheres in $\mathbb{R}^d$. Common choices of
kernels in surrogate optimization are the linear kernel $\varphi(r)=r$,
the cubic kernel $\varphi(r)=r^3$,
and the thin-plate spline $\varphi(r)=r^2\log(r)$.
The coefficients $\lambda_i$ are determined by imposing the
interpolation conditions $s_{f,X}(x_i) = f(x_i)$ for $i=1,\ldots,n$ and
the discrete orthogonality condition
\begin{equation}
  \label{eq:disc_orthg}
  \sum_{i=1}^n \lambda_i q(x_i) = 0, \qquad \forall q \in \Pi_{k-1}^d.
\end{equation}
If we let $\{\pi_i\}_{i=1}^{m}$ be a basis for the
$m=\binom{k-1+d}{k-1}$-dimensional linear space $\Pi_{k-1}^d$, so we
can write $p(x) = \sum_{i=1}^m c_i  \pi_i(x)$, the interpolation
conditions lead to the following linear system of equations
\begin{equation}
    \label{eq:rbf_system}
    \begin{bmatrix} \Phi & P \\ P^T & 0  \end{bmatrix}
    \begin{bmatrix} \lambda \\ c \end{bmatrix}=
    \begin{bmatrix} f_X \\ 0 \end{bmatrix},
\end{equation}
where $\Phi_{ij} = \varphi(\|x_i - x_j\|)$, $P_{ij} = \pi_j(x_i)$, and $f_X = [f(x_1), \ldots, f(x_n)]^T$.
The solution to the linear system of equations is unique as long as $\text{rank}(P) = m$
and $k$ is at least the order of the kernel $\varphi$. The cubic and thin-plate spline
kernels are both of order $k=2$, so a polynomial tail of at least degree 1 is necessary, which is what we use.

A direct solver of the RBF system requires computing a dense LU factorization
at a cost of $\calO(n^3)$ flops. We can utilize the fact that we are adding
a few points at a time, which allows incremental updates of an initial factorization
in quadratic time. We first evaluate $n$ points such that $\text{rank}(P) = m$, which
makes it possible to compute an initial LU factorization with pivoting
\[
    A = \begin{bmatrix} 0 & P^T \\ P & \Phi \end{bmatrix} = PL_{11}U_{11},
\]
where we have reordered the blocks to make it more natural to add new points to
the system. After adding the $k$ new points $\hat{X} = \{\hat{x}_i\}_{i=1}^{k}$
we want to find an LU factorization of the extended system
\[
    \hat{A} = \left[\begin{array}{cc|c}
                  0         & P^T          & \hat{P}^T     \\
                  P         & \Phi         & \hat{\Phi}    \\
                  \hline
                  \hat{P}   & \hat{\Phi}   & \hat{\varphi} \rule{0pt}{2.6ex}\\
                \end{array}\right]
              :=
                \begin{bmatrix} A & B \\ B^T & C \end{bmatrix}
\]
where $\hat{\Phi}_{ij} = \varphi(\|x_i - \hat{x}_j\|)$,
$\hat{P}_{ij} = \pi_j(\hat{x}_i)$, and  $\hat{\varphi}_{ij} = \varphi(\|\hat{x}_i - \hat{x}_j\|)$.
The fact that the
trailing Schur complement is positive semi-definite allows us to look for a factorization
of the form
\[
    \begin{bmatrix} A & B \\ B^T & C \end{bmatrix} =
        \begin{bmatrix} P & 0 \\ 0 & I \end{bmatrix}
        \begin{bmatrix} L_{11} & 0 \\ L_{21} & L_{22} \end{bmatrix}
        \begin{bmatrix} U_{11} & U_{12} \\ 0 & L_{22}^T \end{bmatrix} =
        \begin{bmatrix} PL_{11}U_{11} & PL_{11}U_{12} \\
            L_{21}U_{11} & L_{21}U_{12}+L_{22}L_{22}^T \end{bmatrix}.
\]
We need to solve the two triangular systems $B=PL_{11}U_{12}$ and $B^T=L_{21}U_{11}$
followed by computing a Cholesky factorization of $C-L_{21}U_{12}$. This allows us to update the factorization
in $\mathcal{O}(kn^2)$ flops, which is
better than computing a new LU factorization in $\mathcal{O}(n^3)$ flops.

In practice, we add regularization to the system
by using the kernel $\tilde{\varphi}(x_i,x_j) = \varphi(x_i,x_j) + \eta\delta_{ij}$,
for some regularization parameter $\eta \geq 0$, which ensures that the trailing
Schur complement is positive definite and that the system is well-conditioned.

\subsubsection{Gaussian processes}
\label{sec:gp}
A Gaussian process (GP) is stochastic process where any finite
number of random variables have a joint Gaussian distribution;
see,~e.g.~\cite{rasmussen2006gaussian}.
This defines a
distribution over functions $f(x) \sim \mathcal{GP}(\mu(x),k(x,x'))$,
where $\mu : \bbR^d \rightarrow \bbR$ is the mean function and
$k : \bbR^d \times \bbR^d \rightarrow \bbR$ is the covariance kernel.
The GP model allows predicting the value and variance at any point,
so it gives us an idea about the uncertainty of the prediction.
The most popular kernel is the squared exponential kernel
$k(x,y) = \exp(-0.5\|x-y\|^2/\ell^2)$, and other possibilities include
the Mat\'ern kernels. For any points $X = \{x_1,\ldots,x_n\} \subset \bbR^d$,
$f_X \sim \mathcal{N}(\mu_X, K_{XX})$
where $\mu_X$ denotes the vector values for $\mu$
evaluated at each of the $x_i$, and $(K_{XX})_{ij} =k(x_i, x_j)$.
We assume we observe function values $y_X \in \bbR^n$, where
each entry is contaminated by independent Gaussian noise with constant variance $\sigma^2$.
Under a Gaussian process prior depending on the hyper-parameters $\theta$,
the log marginal likelihood is given by
\begin{equation}
  \label{eq:mloglik}
  \mathcal{L}(y_X \, | \, \theta) =
  -\frac{1}{2}\left[(y-\mu_X)^T\alpha + \log |\tilde{K}_{XX}| + n\log 2\pi\right]
\end{equation}
where $\alpha = \tilde{K}_{XX}^{-1}(y_X - \mu_X)$ and $\tilde{K}_{XX} = K_{XX} + \sigma^2 I$ ($\sigma=0$ for a
deterministic $f(x)$). Optimization of (\ref{eq:mloglik}) is expensive, since direction computation of
$\log |\tilde{K}_{XX}|$ involves computing a Cholesky factorization of $\tilde{K}_{XX}$.
The iteration cost of $\calO(n^3)$ quickly becomes significantly more expensive than
using an RBF interpolant, even though both methods are based on kernel interpolation,
and the dependency of the hyper-parameters stops us from updating a factorization when
adding new points. Promising work on scalable approximate Gaussian process regression
can decrease the kernel learning to $\calO(n)$
\cite{wilson2015kernel,dong2017scalable}, but these ideas only work in low-dimensional spaces.
The computational cost for computing the surrogate
model should be compared to the computational cost of function evaluations as we want to spend
most of the computational effort on doing function evaluations.

\subsubsection{Other choices}
RBFs and GPs are by far the two most popular surrogate models in computationally
expensive optimization. We briefly mention some other possible surrogate models
available in \pysot, even though they are not as frequently used. Multivariate
adaptive regression splines (MARS)
\cite{friedman1991multivariate}, are also weighted sums of basis functions
$B_i(x)$, where each basis function is either constant and equal to 1, a hinge function
of either the form $\max(0,x-c)$ or $\max(0,c-x)$ for some constant $c$, or a
product of hinge functions. It is also possible to use polynomial regression
or support vector regression (SVR).
Multiple surrogate models can be combined into
an ensemble surrogate and Dempster-Shafer theory can be used to decide how to weigh
the different models \cite{muller2011mixture}. This is useful in situations where it is
hard to know what surrogate model to choose for a specific problem. \cite{muller2014influence} indicated regression
polynomial surrogate did not perform well by themselves on test problems, but they were sometimes helpful in combination
with RBF surrogates.

\subsection{Auxiliary problem}
Evaluation of $f(x)$ is expensive, so we optimize an acquisition
function $\alpha(x)$ involving the surrogate model and previously evaluated points
to find the next point(s) to evaluate. We refer to the optimization of $\alpha$ as an auxiliary problem.
This auxiliary problem
must balance exploration and exploitation, where exploration emphasizes evaluating points far from previous
evaluations to improve the surrogate model and escape local minima, while exploitation aims to improve promising
solutions to make sure we make progress. The subsections below describe methods in \pysot to solve the auxiliary problem.

\subsubsection{Candidate points}
An acquisition function based on the weighted-distance merit function is introduced in \cite{regis2007stochastic}
to balance exploration and exploitation. The main idea is to generate a set of
candidate points $\Omega$ and use the merit function
to pick the most promising candidate points. Exploration is achieved by giving preference
to candidate points far from previous evaluations. More specifically,
for each $x\in\Omega$ we let $\Delta(x)$ be the distance
from $x$ to the point closest to $x$ that is currently being or has been evaluated.
By defining $\Delta^{\max} = \max\{\Delta(x) : x \in \Omega\}$ and
$\Delta^{\min} = \min\{\Delta(x) : x \in \Omega\}$ a good measure of exploration
is a small value of $V^D(x) = \frac{\Delta^{\max} - \Delta(x)}{\Delta^{\max}-\Delta^{\min}}$,
where $0 \leq V^D(x) \leq 1$ for all $x \in \Omega$.
Exploitation is achieved through the surrogate model $s(x)$, where a small value
of the quantity $V^S(x) = \frac{s(x) - s^{\min}}{s^{\max}-s^{\min}}$
provides a measure of exploitation, where $s^{\max} = \max\{s(x) : x \in \Omega\}$
and $s^{\min} = \min\{s(x) : x \in \Omega\}$.

The best candidate point is the minimizer of the acquisition function, for a given
$w \in [0,1]$. This shows that $w$ serves as a balance between exploitation and
exploration. A weight close to 0 emphasizes exploration while a weight close
to 1 emphasizes exploitation. Algorithm \ref{alg:cand_points} shows how to select
the most promising candidate point. (see next page)

\begin{algorithm}
    \caption{Candidate point selection}
    \label{alg:cand_points}
    \begin{algorithmic}[1]
        \State Compute $s^{\max} \leftarrow \max\limits_{x \in \Omega} \,\,s(x)$ and
        $s^{\min} \leftarrow \min\limits_{x \in \Omega} \,\,s(x)$

        \For{each $x \in \Omega$}
            \State $V^S(x)\leftarrow
                \begin{cases}
                    \frac{s(x) - s^{\min}}{s^{\max}-s^{\min}} &\text{ if } s^{\max} > s^{\min}  \\
                    1 &\text{  otherwise }
                \end{cases}$
        \EndFor

        \For{each $x \in \Omega$}
            \State $\Delta(x) \leftarrow \min\limits_{y \in \mathcal{A}} d(x,y)$
        \EndFor
        \State Compute $\Delta^{\max} \leftarrow \max\limits_{x \in \Omega} \,\,\Delta(x)$
            and $\Delta^{\min} \leftarrow \min\limits_{x \in \Omega} \,\,\Delta(x)$

        \For{each $x \in \Omega$}
            \State $V^D(x)\leftarrow
                \begin{cases}
                    \frac{\Delta^{\max} - \Delta(x)}{\Delta^{\max}-\Delta^{\min}} &\text{ if }
                        \Delta^{\max}>\Delta^{\min}  \\
                    1 &\text{  otherwise }
                \end{cases}$
        \EndFor
        \Return $\mathop{\mathrm{argmin}}\limits_{x \in \Omega} \,\,\, wV^S(x) + (1-w)V^D(x)$
    \end{algorithmic}
\end{algorithm}

The LMS-RBF method \cite{regis2007stochastic} is useful for low-dimensional optimization
problems. Given a sampling radius $\sigma$, the candidate points
are generated as $\mathcal{N}(0,\sigma^2)$ perturbations along each coordinate
direction from the best solution \cite{regis2007stochastic}. Large values of the sampling radius
will generate candidate points far away from the best solution while
smaller values of the sampling radius will generate candidate points that are close
to the best solution. We defer a description of how the sampling radius is updated to the
next section. If $\sigma$ is smaller than $1$ for an integer variable,
$\sigma=1$ is used to ensure that this variable is also perturbed \cite{muller2013so}.

The DYCORS method \cite{regis2013combining} was developed for high-dimensional problems and the
idea is to start by perturbing most coordinates and perturb fewer
dimensions towards the end of the optimization run \cite{regis2013combining}.
This is achieved by assigning a probability to perturb each dimension.
If $n_0$ points are used in the experimental design and the
evaluation budget is given by $N_{\max}$, each coordinate is perturbed with
probability $p_{\text{perturb}}(n)$
for $n_0 \leq n \leq N_{\max}$. The
probability function used in \pysot is the one introduced in
\cite{regis2013combining}, which is
$p_{\text{perturb}}(n)=\min\left(\frac{20}{d},1\right)\times
\left[1-\frac{\log(n-n_0)}{\log(N_{\max}-n_0)}\right]$.

In \pysot it is also possible to choose candidate points uniformly from $\mathcal{D}$
and use the merit function $wV^S(x) + (1-w)V^D(x)$ to pick the most promising points, which contrasts
to the previous two methods by not making local perturbations around the current best
solution. This helps diversifying the set of evaluated points but resulting in \cite{regis2007stochastic}
for this approach in GMSRBF were not promising.

\subsubsection{Acquisition functions in Bayesian optimization}
Gaussian processes allow us to use acquisition functions that
takes the prediction variance into account. A popular choice is the probability
of improvement, which takes the form
\begin{equation}
 \text{PI}(x) = P(f(x) \leq f(x^{+}) - \xi) = \Phi\left(\frac{f(x^{+}) - \mu(x) - \xi}{\sigma(x)}\right)
\end{equation}
where $\xi$ is a trade-off parameter that balances exploration and exploitation. With $\xi=0$,
probability of improvement does pure exploitation. A common choice is to start with
large $\xi$ and lower $\xi$ towards the end of the optimization run \cite{brochu2010tutorial}.

Expected improvement is likely the most widely used acquisition function in Bayesian optimization, where the
main idea is choosing the point that gives us the largest expected improvement. Mo\u{c}kus
defined improvement as the function
\begin{equation}
    I(x) = \max\{0, f(x^{+}) - f_{n+1}(x)\},
\end{equation}
which can be evaluated analytically under a Gaussian process posterior and
\citet{jones1998efficient} shows that
\begin{equation}
    \text{EI}(x) = \begin{cases}
(f(x^{+}) - \mu(x))\Phi(Z) + \sigma(x)\varphi(Z) & \text{if } \sigma(x) > 0 \\
0                                                & \text{if } \sigma(x) = 0
\end{cases}
\end{equation}

where $Z = (f(x^{+}) - \mu(x))/\sigma(x)$. We can in a similar fashion add a trade-off
parameter $\xi$ in which case the expected improvement takes the form
\begin{equation}
    \text{EI}(x) = \begin{cases}
        (f(x^{+}) - \mu(x) - \xi)\Phi(Z) + \sigma(x)\varphi(Z) & \text{if } \sigma(x) > 0 \\
        0                                                      & \text{if } \sigma(x) = 0
        \end{cases}
\end{equation}
where $Z = (f(x^{+})-\mu(x)-\xi)/\sigma(x)$. Another option that has been proposed is the
lower confidence bound (LCB)
\begin{equation}
  \text{LCB}(x) = \mu(x) - \kappa \sigma(x),
\end{equation}
where $\kappa$ is left to the user.

\subsubsection{Other choices}
Selecting the point that minimizes the bumpiness of a radial basis function, a concept
introduced by \citet{gutmann2001radial}, is supported in other softwares such as
RBFOpt \cite{costa2014rbfopt}. The knowledge gradient acquisition function introduced
used by \citet{frazier2008knowledge} is implemented in Cornell-MOE
\footnote{\url{https://github.com/wujian16/Cornell-MOE}}.

\section{The asynchronous algorithm}
\label{sec:async}
A surrogate optimization in the flavor of Algorithm \ref{alg:sync_algo} is
easy to implement, but may be
inefficient if the evaluation time is not constant. This can be
because evaluating the simulation model requires larger computational
efforts for some input values (e.g., evaluating a PDE-based objective function
may require more time steps for some values of the decision variable $x$).
Computation time can also vary because of variation in the computational resources.
Dealing with potential function evaluation crashes is less obvious in a synchronous framework, 
where we may either try to re-evaluate or exclude the points from the batch. Finally, dynamically
changing the number of resources is much more straightforward in an asynchronous framework,
which we describe next.

Just as in the synchronous parallel case we start by evaluating an
experimental design. These points can be evaluated asynchronously, but the fact
that we want to evaluate all design points before proceeding to the adaptive phase
introduces an undesirable barrier. This becomes an issue if there are straggling workers or
if some points take a long time to evaluate. The most natural solution 
is to generate an experimental design that makes it possible to let
workers proceed to the adaptive phase once all points
are either completed or pending. To be more precise, assume that we have $p$
workers and that $q$ points are needed to construct a surrogate
model. We can then generate an experimental design with $\geq p+q-1$ points,
which will allow workers to proceed to the adaptive phase once there are no outstanding
evaluations in the experimental design. The adaptive phase differs from Algorithm 
\ref{alg:sync_algo} in the sense that we propose a new evaluation as soon as a 
worker becomes available.

We use an event-driven framework where the master drives the event
loop, updates the surrogate, solves the auxiliary problem, etc., and we have $p$
workers available to do function evaluations. The workload of the master
is significantly less than of the workers, so we can use the same number of
workers as we have available resources. The asynchronous algorithm is illustrated
in Algorithm \ref{alg:async_algo}.

\begin{algorithm}
	\caption{Asynchronous surrogate optimization algorithm}
	\label{alg:async_algo}
  	\begin{algorithmic}[1]
	  	\State \textbf{Inputs:} Initial design points, $p$ workers
	  	\State $X \leftarrow \emptyset$, $f_X \leftarrow \emptyset$
		\State Queue initial design
		\State 

		\For{each worker} \Comment{Launch initial evaluations}
			\State Pop point from queue and dispatch to worker
		\EndFor

		\State 
		\While{stopping criterion not met (e.g., $t < t_{\max}$)} \Comment{Event loop}
			\If{evaluation completed}
				\State Update $X$ and $f_X$
			\ElsIf{evaluation failed and retry desired}
				\State Add back to evaluation queue
			\EndIf
			\State 

			\If{worker ready}
				\If{evaluation queue empty}
					\State Build surrogate model
					\State Solve auxiliary problem and add point to the queue
				\EndIf
				\State Pop point from queue and dispatch to worker
			\EndIf
		\EndWhile

		\State \Return $x_{\text{best}}$, $f_{\text{best}}$
  	\end{algorithmic}
\end{algorithm}

\subsection{Updating the sampling radius in Stochastic SRBF}
\label{sec:sampling_radius}
We now elaborate on how to pick the value of the sampling radius $\sigma$
that is used to generate the candidate points used in the LMS-RBF and DYCORS methods. 
We follow the idea in \cite{regis2007stochastic} where counters
$C_{\text{success}}$ and $C_{\text{fail}}$ are used to track the number of
consecutive evaluations with and without significant improvement. If there are too many failures in serial, the algorithm 
restarts. This idea is extended to synchronous parallel in \cite{regis2009parallel} by processing a batch at a time. 
If $C_{\text{success}}$ reaches a tolerance $\calF_{\text{success}}$ the sampling
radius is doubled and $C_{\text{success}}$ is set to 0. Similarly, if
$C_{\text{fail}}$ reaches  $\calF_{\text{fail}}$ the sampling
radius is halved and $C_{\text{fail}}$ is set to 0.

In the asynchronous setting, we update the counters after each completed function 
evaluation. We do not update the
counters for evaluations that were launched before the last time the sampling
radius was changed. The reason for this is that these evaluations are based on
outdated information. The logic for updating the sampling radius and the best
solution can be seen in Algorithm \ref{alg:async_algo_radius}.

\begin{algorithm}[!t]
	\caption{Sampling radius adjustment routine}
	\label{alg:async_algo_radius}
	\begin{algorithmic}[1]
		\State \textbf{Inputs:} $\sigma$, $f(x_i)$, $x_i$, $f_{\text{best}}$,
								$x_{\text{best}}$, $C_{\text{success}}$,
								$C_{\text{fail}}$, $\calF_{\text{success}}$,
								$\calF_{\text{fail}}, \delta$
			\If{$f(x_i) < f_{\text{best}}$}
				\State $f_{\text{best}} \leftarrow f(x_i)$
				\State $x_{\text{best}} \leftarrow x_i$
				\If{$f(x_i) < f_{\text{best}} - \delta |f_{\text{best}}|$}
					\State $C_{\text{succ}} \leftarrow C_{\text{succ}} + 1$
					\State $C_{\text{fail}} \leftarrow 0$
				\EndIf
			\Else
				\State $C_{\text{succ}} \leftarrow 0$
				\State $C_{\text{fail}} \leftarrow C_{\text{fail}} + 1$
			\EndIf
			\State
			\If{$C_{\text{succ}} = \calF_{\text{succ}}$ \OR
				  $C_{\text{fail}} = \calF_{\text{fail}}$}
				\State $C_{\text{succ}} \leftarrow 0$
				\State $C_{\text{fail}} \leftarrow 0$
				\If{$C_{\text{succ}} = \calF_{\text{succ}}$}
					\State $\sigma \leftarrow \min(2\sigma, \sigma_{\max})$
				\Else
					\State $\sigma \leftarrow \max(\sigma/2, \sigma_{\min})$
				\EndIf
			\EndIf
			\State
			\State \Return $\sigma$, $f_{\text{best}}$, $x_{\text{best}}$
						   $C_{\text{success}}$, $C_{\text{fail}}$
  	\end{algorithmic}
\end{algorithm}

We also follow the recommendations in \cite{regis2007stochastic} and
\cite{regis2009parallel} to restart the algorithm when we reach a maximum
failure tolerance parameter $\calM_{\text{fail}}$ or when the sampling radius $\sigma$
drops below $\sigma_{\min}$. Restarting has shown to be successful for LMS-RBF and DYCORS
as it can be hard to make progress when the surrogate is very biased towards the
current best solution and we may be stuck in a local minimum that is hard to escape. 
Restarting the algorithm can help avoid this issue. We do not terminate pending evaluations
after a restart occurs, but they are not incorporated in the surrogate model or used to
adjust the sampling radius when they finish.

\section{POAP implementation}
\label{sec:poap}
This section describes the Plumbing for Optimization with Asynchronous 
Parallelism\footnote{\poap can be downloaded from the GitHub repository:
https://github.com/dbindel/\poap} (\poap) framework. \poap has three main components: 
a controller that asks workers to run function evaluations, a
strategy that proposes new actions, and a set of workers that carry out
function evaluations.

\subsection{Controllers}
The controller is responsible for accepting or rejecting proposals by the strategy
object, controlling and monitoring the workers, and informing the strategy object
of relevant events. Examples of relevant events are the processing of a proposal,
or status updates on a function evaluation. Interactions between controller and
the strategies are organized around proposals and evaluation records. At the
beginning of the optimization and on any later change to the system state, the
controller requests a proposal from the strategy. The proposal consists of an
action (evaluate a function, kill a function, or terminate the optimization), a
list of parameters, and a list of callback functions to be executed once the
proposal is processed. The controller then either accepts the proposal (and sends
a command to the worker), or rejects the proposal.

When the controller accepts a proposal to start a function evaluation, it creates
an evaluation record to share information about the status of the evaluation with
the strategy. The evaluation record includes the evaluation point, the status of
the evaluation, the value (if completed), and a list of callback functions to be
executed on any update. Once a proposal has been accepted or rejected, the
controller processes any pending system events (e.g. completed or canceled
function evaluations), notifies the strategy about updates, and requests the next
proposed action.

\poap comes with a serial controller for when
objective function evaluations are carried out in serial. There is also a
threaded controller that dispatches work to a set of workers where each worker
is able to handle evaluation and kill requests. The requests are asynchronous
in the sense that the workers are not required to complete the evaluation or
termination requests. The worker is forced to respond to evaluation requests, but may
ignore kill requests. When receiving an evaluation request, the worker should either
attempt the evaluation or mark the record as killed. The worker sends status
updates back to the controller by updating the relevant record. There is also
an extension of the threaded controller that works with MPI and a controller 
that uses simulated time. The latter is useful for testing 
asynchronous optimization strategies for different evaluation time distributions.

\subsection{Strategies}
The strategy is responsible
for choosing new evaluations, killing evaluations, and terminating the optimization
run when a stopping criteria is reached. \poap provides some basic default strategies
based on non-adaptive sampling and serial optimization routines and also some
strategies that adapt or combine other strategies.

Different strategies can be composed by combining their control actions, which
can be used to let a strategy cycle through a list of optimization strategies
and select the most promising of their proposals. Strategies can also subscribe
to be informed of all new function evaluations so they incorporate any new function
information, even though the evaluation was proposed by another strategy.
This makes it possible to start several independent strategies while still
allowing each strategy to look at the function information that comes from
function evaluations proposed by other strategies. As an example, we can have a
local optimizer strategy running a gradient based method where the starting
point can be selected based on the best point found by any other strategy.
The flexibility of the \poap framework makes combined strategies like these
straightforward.

\subsection{Workers}
The multi-threaded controller employs a set of workers that are capable of
managing concurrent function evaluations. Each worker does not provide
parallelism on its own, but the worker itself is allowed to exploit parallelism
by separate external processes.

The basic worker class can call Python objective functions, which only results 
in parallelism if the objective
function itself allows parallelism. There is also a worker class that uses
subprocesses to evaluate objective functions that are
not necessarily in Python. The user is responsible for specifying how to evaluate
the objective function and parse partial information.

The number of workers can be adjusted dynamically during the optimization
process, which is particularly useful in a cloud setting. \poap supports running both
on the Google Cloud platform (GCP) and the Amazon Web Services (AWS). We support 
workers connecting to a specified TCP/IP port to communicate with the controller,
making it easy to add external resources.

\section{\lowercase{py}SOT implementation}
\label{sec:pysot}
\begin{figure}[!h]
	\centering
    \includegraphics[width=0.7\textwidth]{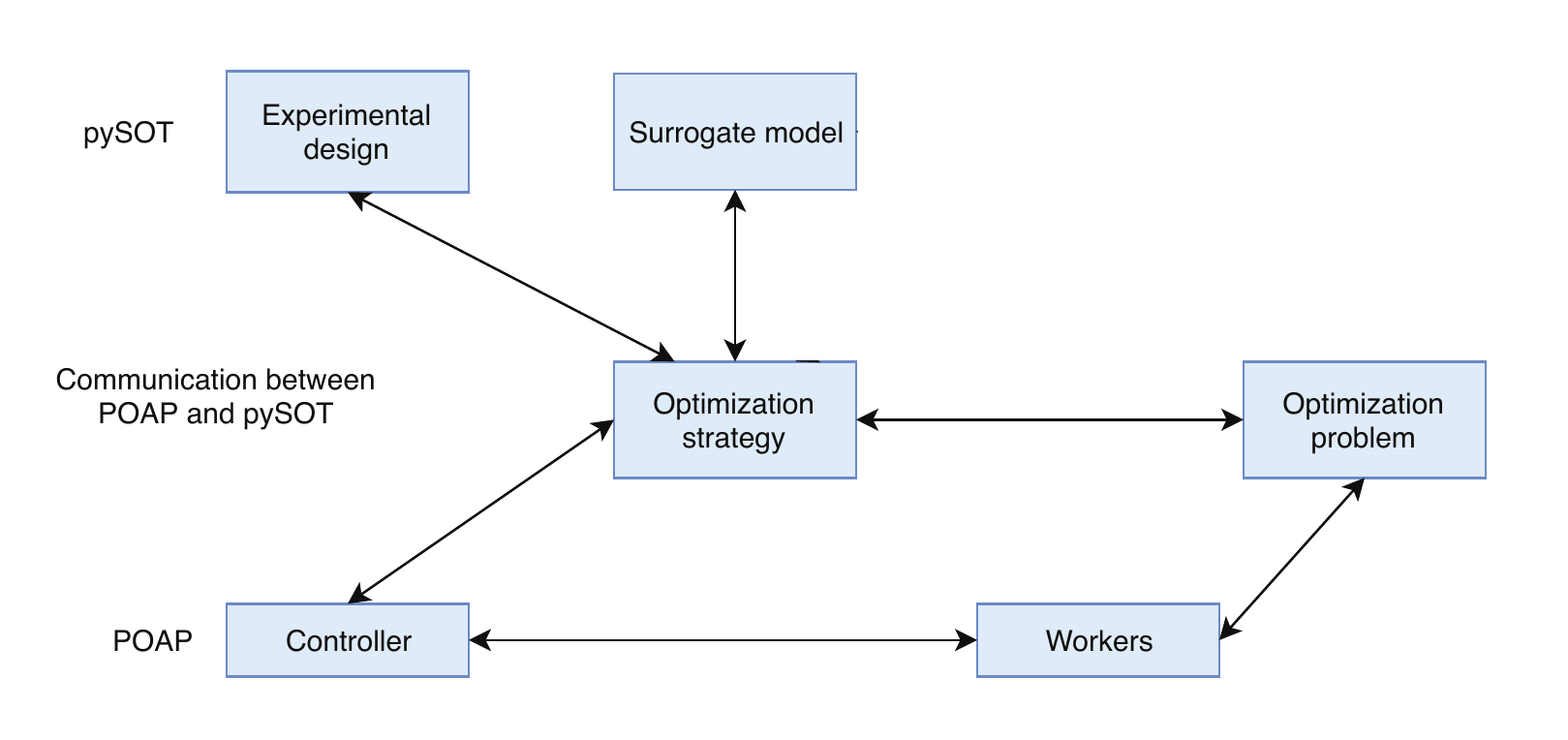}
    \caption{Communication between \poap and \pysot}
    \label{fig:flowchart}
\end{figure}

The surrogate optimization toolbox\footnote{\pysot can be downloaded from the
Github repository: https://github.com/dme65/\pysot} (\pysot) is a collection of
surrogate optimization strategies that can be used with the \poap framework. \pysot
follows the general surrogate optimization framework in Algorithm \ref{alg:sync_algo} 
and allows using asynchrony as was described in Algorithm \ref{alg:async_algo}. 
We illustrate the communication between \poap and \pysot in Figure \ref{fig:flowchart}.
\pysot communicates with \poap through the optimization strategy where the \pysot strategy
is responsible for proposing an action when different events happen. All of the 
worker communication is handled by the \poap controller.

The \pysot objects follow an abstract class definition to make sure that
custom implementations fit the framework design. This is achieved by forcing the
objects to inherit an ABC object design, which makes it easy for users to
add their own implementations. We proceed to describe each object and their 
role in the \pysot framework. 

\subsection{Strategies and auxiliary problems}
The strategy object follows the \poap framework. \pysot
implements an asynchronous base class for surrogate optimization which 
serves as a template for all surrogate optimizations in \pysot.
This base class abstracts out the difference between serial, synchronous parallel, 
and asynchronous parallel. \pysot supports the candidate point methods 
SRBF and DYCORS. We also support strategies for the most common 
acquisition functions from BO: expected improvement (EI) and the lower confidence 
bound (LCB). 

\subsection{Experimental design}
\pysot implements the symmetric Latin hypercube (SLHD), Latin hypercubes (LHD), and 2-factorial
designs that were described in \S\ref{sec:expdes}.
The experimental design is always evaluated first and the asynchronous optimization
strategy in \pysot is designed to proceed to the adaptive phase as soon as no initial
design points are outstanding. Another possibility is to cancel the pending evaluations
from the initial phase and proceed to the adaptive phase as soon as possible, but we choose
to finish the entire initial design as exploration is important for multi-modal optimization
problems. As discussed in the previous section, we must choose
enough initial design point to allow building the surrogate model when all points in the initial
design are either completed or pending.

\subsection{Surrogate models}
\pysot supports the many popular surrogate models, including RBFs, GPs, MARS, polynomial regression, 
and support vector regression. We provide our own RBF implementation that uses the incremental
factorization update idea that was described in \S\ref{sec:rbf}. 
Support for MARS is provided via py-earth\footnote{\url{https://github.com/scikit-learn-contrib/py-earth}} 
and support for GPs and polynomial regression is 
provided through scikit-learn \cite{pedregosa2011scikit}. The surrogate model does not need
access to any of the other objects, as it just constructs a model based on the evaluated points
and their values. The surrogate fitting problem may be ill-conditioned if the domain is scaled poorly,
and we provide wrappers for rescaling the domain to the unit hypercube, which is particularly
useful on problems where the bounds are very skewed. We add regularization to the linear system 
when radial basis functions are used to keep the system well-conditioned. Previous work has shown
that hard-capping of function values can be useful to avoid oscillation, where a common choice is to 
replace all function values above the median by the median function value, and we provide
wrappers for this as well.

\subsection{Optimization problems}
The optimization problem object specifies the number of dimensions, the number of analytical 
constraints, and provide methods for evaluating the objective function and the constraints. 
We provide implementations of many standard test problems which can be used to compare 
algorithms within the \pysot framework. The optimization problem does not depend on any 
other objects.

\subsection{Checkpointing}
Checkpointing is important when optimizing an expensive function since the optimization
may run for several days or weeks, and it would be devastating if all information was 
lost due to e.g., a system or power failure. \pysot supports a controller wrapper
for saving the state of the system each time something changes, making it possible to resume 
from the latest such snapshot.

\section{Code examples}
\label{sec:code}
This section illustrates how \poap and \pysot can be used to minimize the
Ackley test function. We will use the threaded controller and asynchronous
function evaluations. The code is based on \pysot version 0.2.3 and
\poap version 0.1.26.

Our goal in this example is to minimize
the 10-dimensional Ackley function, which is a common test function in global optimization.
We use a symmetric Latin hypercube, an RBF surrogate with a cubic kernel and linear tail, and
the DYCORS strategy for generating candidate points. Importing the necessary modules can be
done as follows:
\begin{python}
from pySOT.optimization_problems import Ackley
from pySOT.experimental_design import SymmetricLatinHypercube
from pySOT.surrogate import RBFInterpolant, CubicKernel, LinearTail
from pySOT.strategy import DYCORSStrategy
from poap.controller import ThreadController, BasicWorkerThread
\end{python}
We next create objects for the optimization problem, experimental design,
and surrogate model.
\begin{python}
num_threads = 4
max_evals = 500

dim = 10
ackley = Ackley(dim=dim)
rbf = RBFInterpolant(
        dim=dim, kernel=CubicKernel(), tail=LinearTail(dim))
slhd = SymmetricLatinHypercube(dim=dim, num_pts=2*(dim+1))
\end{python}
We are now ready to launch a threaded controller that will run
asynchronous evaluations. We create an instance of the DYCORS strategy and
append it to the controller.
\begin{python}
controller = ThreadController()
controller.strategy = DYCORSStrategy(
        opt_prob=ackley, exp_design=slhd, surrogate=rbf,
        max_evals=max_evals, asynchronous=True)
\end{python}
We need to launch the workers that do function evaluations.
In this example we use standard threads and give each worker an
objective function handle.
\begin{python}
for _ in range(num_threads):
    worker = BasicWorkerThread(controller, ackley.eval)
    controller.launch_worker(worker)
\end{python}
The workers have been launched and the optimization strategy has been created,
so we are ready to start the optimization run. The following code runs the
optimizer and prints the best solution.
\begin{python}
result = controller.run()
print("Best value found: {0}".format(result.value))
print("Best solution found: {0}".format(result.params[0]))
\end{python}

\section{Numerical experiments}
\label{sec:experiments}
In this section we study the performance of serial, synchronous parallel,
and asynchronous parallel when varying the evaluation time distribution and the
number of processors. We focus on the DYCORS method using a cubic RBF
surrogate with a linear tail as both have very low overhead, allowing us to run
many trials, each with a large number of function evaluations. Previous work
has shown that DYCORS outperforms the competing methods for computationally expensive
multi-modal functions in a large number of dimensions \cite{regis2013combining}.

The evaluation times are drawn from a Pareto distribution
with probability density function (PDF) given by:
\begin{equation}
  \frac{\alpha b^{\alpha}}{x^{\alpha+1}}\chi_{[b,\infty)}(x).
\end{equation}
The Pareto distribution is heavy-tailed for small values of $\alpha$ and this case
is suitable for studying large variance in the evaluation time.
We use $b = 1$ so the support is $[1,\infty)$ and use
different values of $\alpha$ to achieve different tail behaviors. This setup
models homogeneous resources and spatial dependence.
We use $\alpha \in \{102, \,12,\, 2.84\}$ which corresponds to standard deviations
$0.01$, $0.1$, and $1$.

We run the serial and synchronous parallel versions with their default
hyper-parameters since both methods showed good results in \cite{regis2007stochastic}
and \cite{regis2009parallel} respectively. The hyper-parameter values used for the asynchronous
algorithm are the same as for the synchronous parallel version except for
$\calF_{\text{fail}}$ which we multiply by $p$ since we count evaluations rather
than batches. The hyper-parameter values are shown in Table \ref{tab:paramvals}.
We restart the algorithm with a new experimental
design if at some point $\sigma = \sigma_{\min}$ and the algorithm has
failed to make a significant improvement in the last $\calM_{\text{fail}}$
evaluations. Our experiments use $4$, $8$, $16$, and $32$ workers for the parallel
algorithms and we give each algorithm an evaluation budget of
$50*32=1600$ evaluations. This is an upper bound for a time budget of $50$ units
of time. We exclude the overhead
from fitting the surrogate and generating candidate points since this is
negligible when the function evaluations are truly expensive.

\begin{table*}[b]
    \centering
    \resizebox{0.8\textwidth}{!}{
    \begin{tabular}{lc}
		Hyper-parameter & Value \vspace{-0.5mm} \\
		\hline \vspace{-3mm} \\
		$|\Omega_n|$ (number of candidate points per proposal) & $100d$ \vspace{0mm} \\
		$\Upsilon$ (weight pattern) & $\langle 0.3, 0.5, 0.8, 0.95 \rangle$ \vspace{0mm} \\
		$\kappa$ (number of weights in $\Upsilon$) & $4$ \vspace{0mm} \\
		$\sigma_{\text{init}}$ (Initial step size) & $0.1 \ell (\calD)$ \vspace{0mm} \\
		$\sigma_{\min}$ (minimum step size) & $0.1(1/2)^6 \ell(\calD)$ \vspace{0mm} \\
		$\delta_{\text{tol}}$ (radius tolerance) & $0.0025 \ell(\calD)$ \vspace{0mm} \\
		$\calF_{\text{succ}}$ (threshold parameter for increasing the step size) & $3$ \vspace{0mm} \\
		$\calF_{\text{fail}}$ (tolerance parameter for decreasing step size) &
			$p \lceil \max(4/p, d/p)\rceil$ \vspace{0mm} \\
		$\calM_{\text{fail}}$ (maximum failure tolerance parameter) & $4\, \calF_{\text{fail}}$ \\
		\hline
	\end{tabular}}
	\caption{Hyper-parameter values used for the asynchronous algorithm}
	\label{tab:paramvals}
\end{table*}

We consider the multimodal test problems F15-F24 from the BBOB test suite
\cite{hansen2009real}. These problems are challenging and non-separable
and we use the $10$-dimensional versions for our experiments. The domain for each
problem is $[-5, 5]^{10}$, and location of the global minimum and the value at
the global optimum are generated randomly depending on what instance is being
used, where we use instance $1$ for each problem.
These problems are not expensive to evaluate, but we pretend they
are computationally expensive black-box functions and draw the evaluation
time from a Pareto distribution.

\begin{figure}[!ht]
	\centering
    \resizebox{0.96\textwidth}{!}{
    	\subfigure{\label{fig:F15}
            \includegraphics[width=0.96\textwidth,trim={0cm 0cm 0cm 0cm},clip]{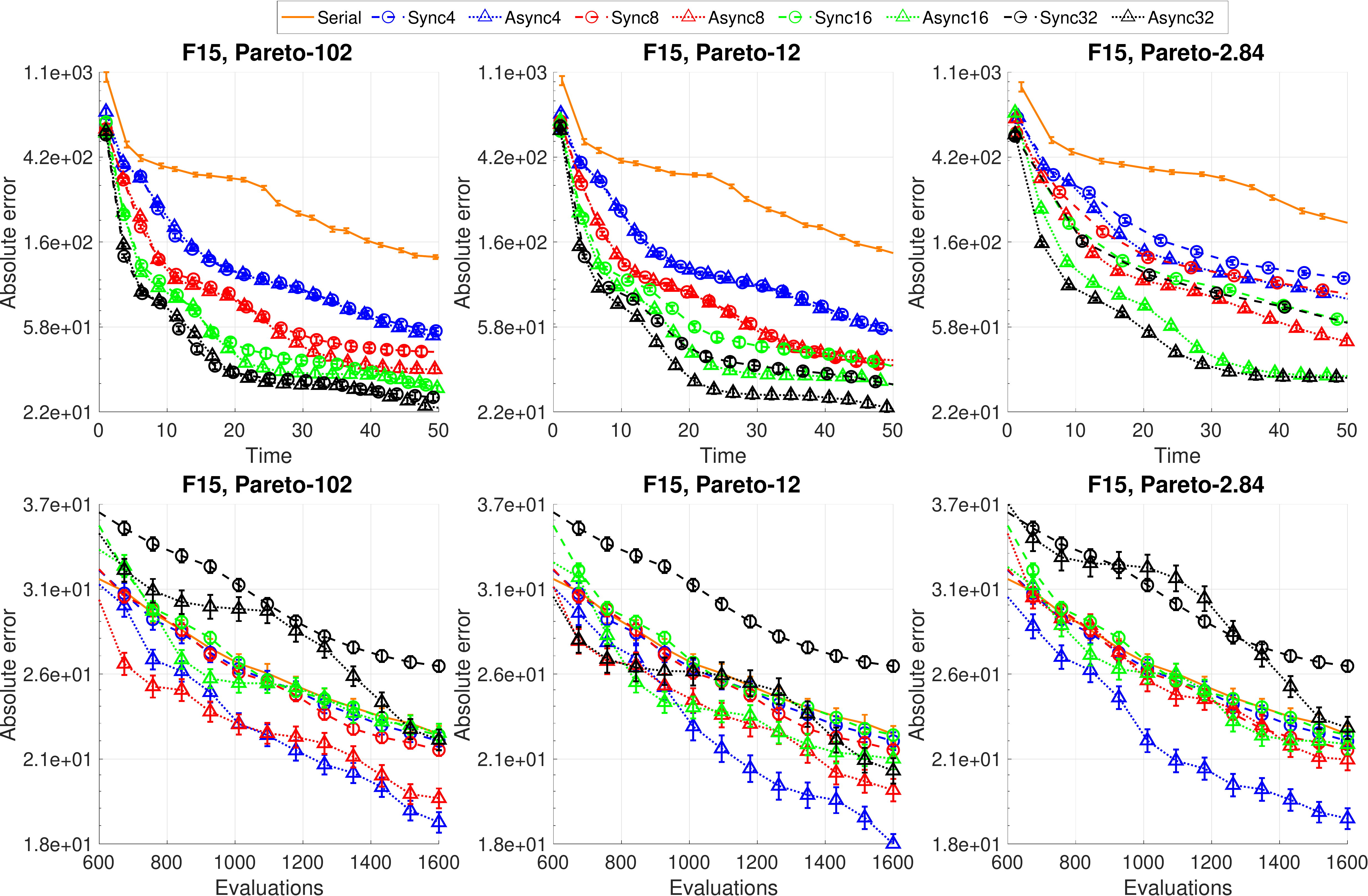}}
    }
    \resizebox{0.96\textwidth}{!}{
        \subfigure{\label{fig:F17}
            \includegraphics[width=0.96\textwidth,trim={0cm 0cm 0cm 0cm},clip]{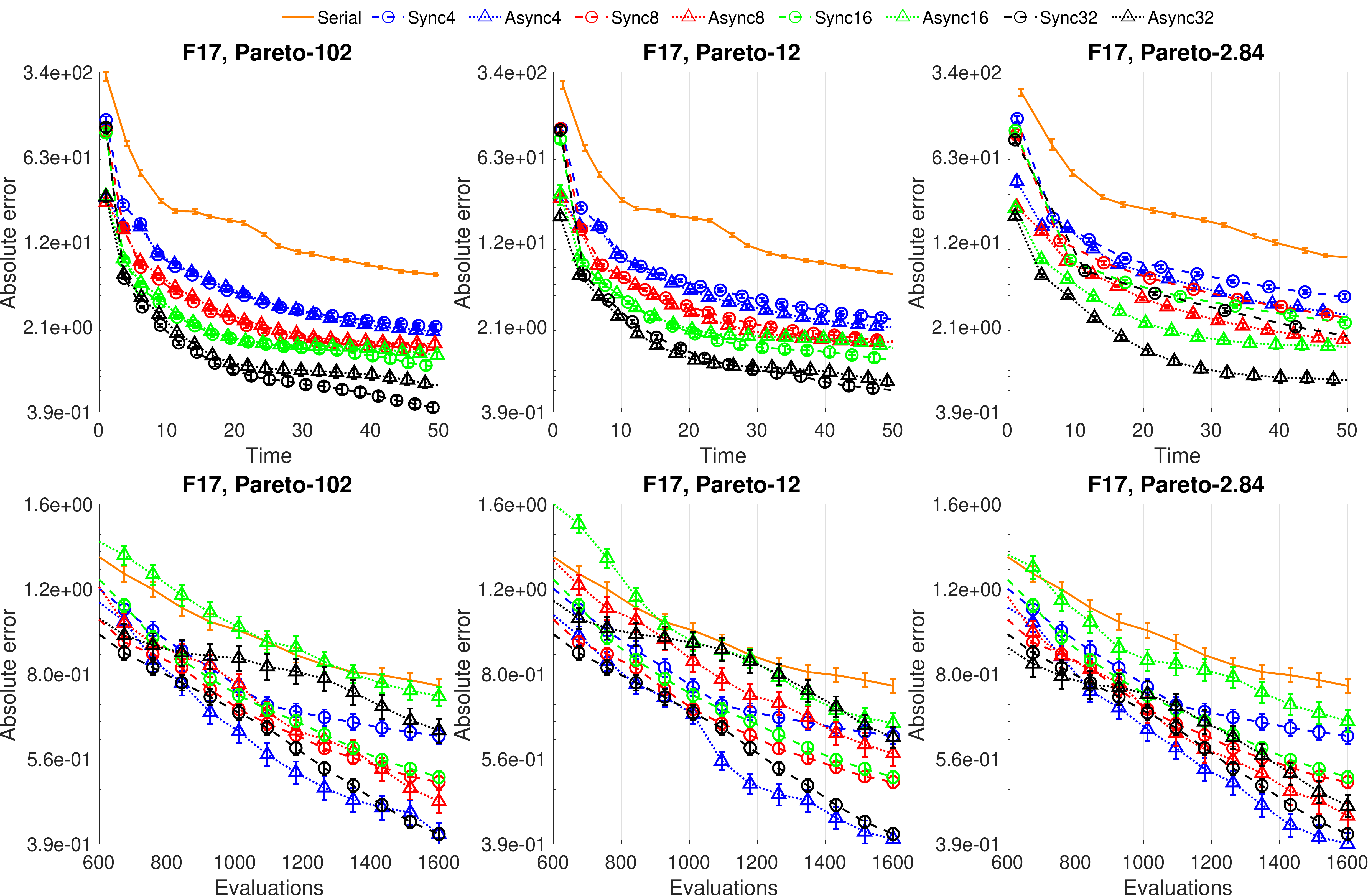}}
    }
    \caption{Progress vs time and progress vs number of evaluation for F15 and F17.
    The error bars show the standard error based on 100 trials using 1600 evaluations.
    The plots with respect to number of evaluations are zoomed in to make the lines easier
    to distinguish.}
    \label{fig:progplot1}
\end{figure}

We compare progress vs time and progress vs number of
evaluations. Comparing progress vs time will show what method does
well in practice since we are often constrained by a time
budget rather than an evaluation budget. We will also be able to see
the effect of adding more processors, which is expected to be fruitful since
exploration is important for multi-modal problems. We compare progress vs number of
evaluations to study the importance of information. The
serial and synchronous methods are independent of the evaluation time in this
case since there is a barrier after each batch. The asynchronous algorithm is
affected by the variance, which affects how much information is available
at a given iteration.
The serial version
always has more points incorporated in the surrogate at a given iteration, but
explores less than the parallel versions.
Figure \ref{fig:progplot1} shows the experimental results for F15 and F17.

These problems are chosen because they show the key points we are
trying to convey. The first row in each plot shows absolute error vs
time and the second row shows absolute error vs number of evaluations.
The error is the difference between the best objective function found so
far and the value of the true optimization. The
absolute error is plotted on a log-scale in all plots to make it easier to
interpret the results. This should be taken into account when looking at the
error bars, which show the standard error of the estimated mean based on
100 trials. Note that each row has the same range in absolute error
to make the comparison easier.

F15 illustrates a case where synchronous and asynchronous parallel perform similarly
with respect to time when the variance is small. The difference grows when
the variance increases and asynchrony is always superior in the large variance
case. Synchrony does slightly better than asynchrony in the small variance case
for F17. The results versus number evaluations are interesting and asynchrony
with 4 processors is consistently the best choice for both F15 and F17. This is
unexpected since the asynchronous versions never have more information than the
corresponding synchronous versions, indicating that maximizing information is
not as important on multi-modal problems. We also see that the serial version
is outperformed by the parallel versions when looking at number of evaluations.
This is another indicator that exploration is more important than information.

We can also compare the speedup from using more workers for
synchronous and asynchronous parallel. Speedup is measured by the quantity
\begin{equation*}
  S(p) = \frac{T^*(1)}{T(p)} = \frac
        {\text{Execution time for fastest serial algorithm to reach target value}}
        {\text{Execution time for parallel algorithm with }p\text{ processors to reach target value}}.
\end{equation*}
This requires knowledge of the fastest serial algorithm which is hard to know
given randomized initial conditions. We therefore consider relative speedup
where $T^*(1)$ is replaced by $T(1)$. We will measure speedup by running
synchronous and asynchronous parallel and compare the results to the serial
case. We need to estimate the expected value of the speedup since all of our
algorithms are stochastic.

\begin{figure}[!ht]
	\centering
    \resizebox{0.9\textwidth}{!}{
    	\subfigure{\label{fig:F15_speedup}
            \includegraphics[width=0.9\textwidth,trim={5cm 0cm 3cm 0cm},clip]{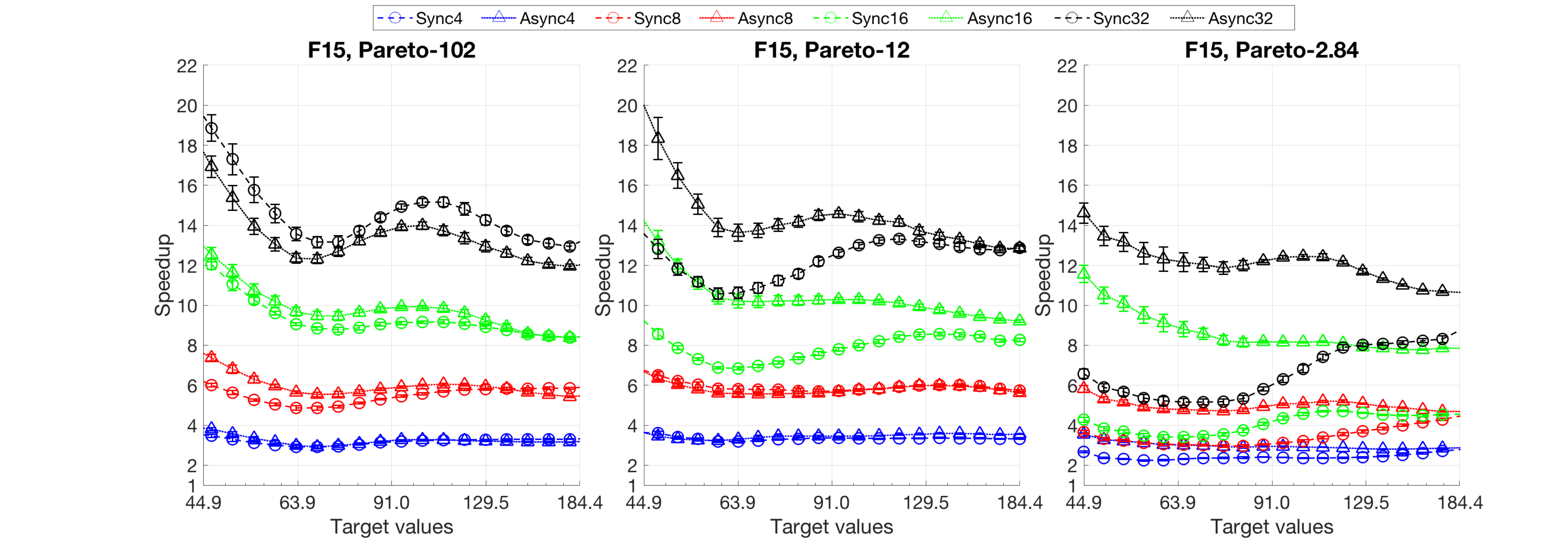}}
    }
    \resizebox{0.9\textwidth}{!}{
        \subfigure{\label{fig:F17_speedup}
            \includegraphics[width=0.9\textwidth,trim={5cm 0cm 3cm 0cm},clip]{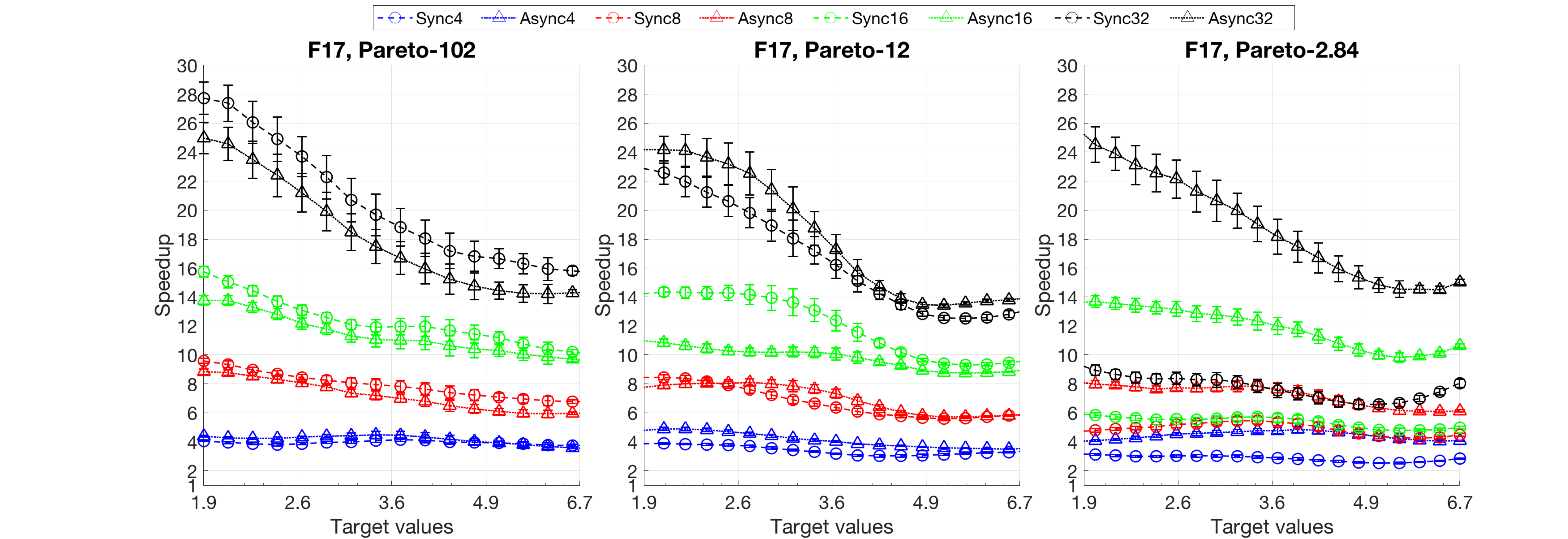}}
    }
    \caption{Relative speedup for different target values for F15 and F17. The error bars
    show the standard error based on 100 trials using 1600 evaluations. Results show the speedup of reaching the target values as in equation (12). }
    \label{fig:speedupplot1}
\end{figure}

A main problem with speedup tests is choosing a good target value
where the speedup is measured. For unimodal problems, such a target can be based
on a small neighborhood of the global minimum, but this is unreasonable for
multimodal test problems. Our approach is to run each algorithm
for 100 trials, compute the intersection of the ranges of function values,
and compute the speedup for a set of targets within the intersection. This
allows us to see how the speedup depends on different target values. Figure
\ref{fig:speedupplot1} shows the speedup for F15 and F17.

We achieve close to linear speedup with asynchrony for small target values when
using 4, 8, and 16 processors on F15 in the case of small variance. The speedup
is larger for synchrony when we consider 32 processors, but is clearly sub-linear.
The speedup is larger for small target values, which indicates that the serial
algorithm is more likely to get stuck in a local minimum which triggers a restart.
The effect of increasing the variance clearly degrades the performance of
synchrony while the results for asynchrony do not change much. The speedup on
F17 is generally better for the synchronous algorithm in the case of small
variance.


\section{Conclusions}
\label{sec:conclusions}
We have introduced the event-driven optimization framework \poap, which provides
an easy way to build and combine new optimization algorithms for computationally
expensive functions. \poap has
three main components, a controller, a strategy, and a collection of workers. The
controller accepts or rejects proposals from the strategies, monitors the workers,
and informs the strategy when new events occur. The strategy proposes actions when
an event occurs, such as starting a new function evaluation, re-evaluating an input,
or terminating the optimization run. The workers do 
function evaluations when instructed by the controller and they support
partial updates, making it possible to terminate
unpromising evaluations. The flexibility of the \poap framework makes it easy to
combine optimization strategies.

We have also introduced the inter-operable library \pysot that supports the most 
popular surrogate optimization methods. \pysot is a
collection of synchronous/asynchronous strategies, experimental designs, surrogate
models, and auxiliary functions, that are commonly used in surrogate
optimization. \pysot comes with a large set of standard test problems and 
efficiently serves as a test
suite in which new optimization algorithms can be compared to existing
methods using asynchronous or synchronous parallel.
The object oriented design makes it easy to add new functionality and there
is also support for resuming crashed and terminated runs.

We have introduced a general asynchronous surrogate optimization method
for computationally expensive black-box problems that extends the work in
\cite{regis2009parallel} to the case when function evaluation times are not
necessarily constant and computational resources are
not necessarily homogeneous. Our version also handles worker failures and
evaluation crashes, which was not considered in \cite{regis2009parallel}.

The numerical experiments show that asynchrony performs similarly to synchrony
even when the variance in evaluation time is small. Comparing progress vs
number of evaluations showed that the serial method, which maximizes the
information at each step, does not outperform the parallel methods. This is
likely because exploring is more important than maximizing the information for
each sample. The implication is that idle time is more important than
information and the asynchronous method clearly outperforms the synchronous method
when we increase the variance in evaluation time or number of processors. We
also studied a unimodal function, in which case the serial method
performs best when we compare progress vs number of evaluations. This is expected
since exploration is fruitless for unimodal functions.

A relative speedup analysis shows good results for the asynchronous method and
we achieve close to linear speedup for 4, 8, and 16 processors. The speedup
for the synchronous method clearly decreases when the variance of the
evaluation time increases, which is expected since this increases idle time.
We conclude that for multi-modal problems asynchrony should be preferred over
synchrony and that adding more processors leads to faster convergence.

\section*{Acknowledgement}
The authors appreciated support from NSF CISE 1116298 to Prof. Shoemaker and Bindel, and Prof.
Shoemaker's start up grant from National University of Singapore. We also thank Dr. Taimoor Akhtar for
his incorporation of multi-objective code into \pysot.

\bibliography{references}
\bibliographystyle{unsrtnat}

\end{document}